%ritsuf.tex: Automatic Counting of Tilings of Skinny Plane Regions 
%%a Plain TeX file by Shalosh B. Ekhad and Doron Zeilberger (x pages)

%begin macros

\baselineskip=14pt
\parskip=10pt

\font\eightrm=cmr8 
\font\eighttt=cmtt8
\magnification=\magstephalf

\def\1{{\overline{1}}}
\def\2{{\overline{2}}}
\parindent=0pt
\overfullrule=0in

\def\frac#1#2{{#1 \over #2}}
%\headline={\rm  \ifodd\pageno  \RightHead  \else  \LeftHead  \fi}
%\def\RightHead{\centerline{
%Title
%}}
%\def\LeftHead{ \centerline{Doron Zeilberger}}
%end macros
\bf
\centerline
{
Automatic Counting of Tilings of Skinny Plane Regions 
}
\rm
\bigskip
\centerline{ {\it
Shalosh B. EKHAD and
Doron 
ZEILBERGER}\footnote{$^1$}
{\eightrm  \raggedright
Department of Mathematics, Rutgers University (New Brunswick),
Hill Center-Busch Campus, 110 Frelinghuysen Rd., Piscataway,
NJ 08854-8019, USA.
%\break
{\eighttt zeilberg  at math dot rutgers dot edu} ,
\hfill \break
{\eighttt http://www.math.rutgers.edu/\~{}zeilberg/} .
Written: June 21, 2012. 
Supported in part by the USA National Science Foundation.
}
}

{\bf Very Important} 

This article is accompanied by the Maple packages

$\bullet$
{\eighttt http://www.math.rutgers.edu/\~{}zeilberg/tokhniot/RITSUF} \quad ,

$\bullet$
{\eighttt http://www.math.rutgers.edu/\~{}zeilberg/tokhniot/RITSUFwt} \quad , 

$\bullet$
{\eighttt http://www.math.rutgers.edu/\~{}zeilberg/tokhniot/ARGF} \quad ,

to be described below. In fact, more accurately, this article
accompanies these packages, written by DZ and the many output files,
discovering and proving deep enumeration theorems, done by SBE, that
are linked to from the webpage of this article

{\tt http://www.math.rutgers.edu/\~{}zeilberg/mamarim/mamarimhtml/ritsuf.html} \quad \quad .

{\bf How It All Started: April 5, 2012}

During one of the Rutgers University Experimental Mathematics Seminar dinners,
the name of Don Knuth came up, and two of the participants, David Nacin,
who was on sabbatical from William Patterson University, 
and first-year graduate student Patrick Devlin,
mentioned that they recently solved a problem that Knuth proposed in {\it Mathematics Magazine}[Kn].
The problem was:

{\bf 1868.} {\it Proposed by Donald E. Knuth, Stanford University, Stanford, CA.}

Let $n \geq 2$ be an integer. Remove the central $(n-2)^2$ squares from an
$(n+2) \times (n+2)$ array of squares. In how many ways can the remaining squares
be covered with $4n$ dominoes?

As remarked in the published solution in Math. Magazine, 
the problem was already
solved in the literature by Roberto Tauraso[T]. The answer turned out to be
very elegant: $4(2F_n^2 +(-1)^n)^2$.

David and Patrick, as well as the 
solution published  in {\it Math. Magazine},
used human ingenuity and {\it mathematical deduction}.
But, as already preached in [Z1] and [Z2], the
following is a {\bf fully rigorous proof}:

``By direct counting of tilings, the 
first $16$ terms (starting at $n=2$) of the enumerating sequence are
$$
36, 196, 1444, 9604, 66564, 454276, 3118756, 21362884, 146458404, 1003749124, 6880038916, 
$$
$$
47155859716, 323212716324, 2215328606404, 15184099435684, 104073336269956 \quad .
$$
But so are the first $16$ terms of the sequence $\{4(2F_n^2 +(-1)^n)^2\}$. Since the statement is true for
the first $16$ terms, it must be true for {\bf all} $n \geq 2$. {\bf QED}!''

In order to {\it justify} this ``empirical'' proof, all we need to say 
is that both sides
are obviously {\it C-finite} sequences  whose minimal recurrences 
have order $\leq 8$, hence
their difference is a C-finite sequence of order $\leq 16$ and hence if it vanishes for
the first $16$ terms, 
it {\it always} vanishes. 
But first let's remind ourselves what are $C$-finite sequences.

{\bf $C$-finite Sequences} 

Recall that a $C$-finite sequence $\{a(n)\}_{n=0}^{\infty}$ is a sequence that satisfies a
{\bf homogeneous linear-recurrence equation with  constant coefficients}. It is known (but not as well-known as it should be!)
and easy to see (e.g. [Z2],[KP]) that the set of $C$-finite sequences is an {\it algebra}.
Even though a $C$-finite sequence is an ``infinite'' sequence, it is in fact, {\bf like everything else in mathematics}
(and elsewhere!) a {\bf finite} object. 
An order-$L$ $C$-finite sequence $\{a(n)\}_{n=0}^{\infty}$
is completely specified by the coefficients $c_1,c_2, \dots , c_L$  of the recurrence
$$
a(n)=c_1 a(n-1)+ c_2 a(n-2) + \dots + c_L a(n-L) \quad,
$$
and the {\bf initial conditions}
$$
a(0)=d_1 \quad , \quad \dots \quad , \quad a(L-1)=d_L \quad.
$$
So a $C$-finite sequence can be {\bf coded} in terms of the $2L$ ``bits'' of information
$$
[[d_1, \dots , d_L],[c_1, \dots, c_L] ] \quad .
$$
For example,  the Fibonacci sequence is written:
$$
[[0,1],[1,1]] .
$$
Since this ansatz (see [Z1]) is fully decidable, it is
possible to decide equality, and  evaluate {\it ab initio}, wide classes of sums, and
things are easier than the {\it holonomic ansatz} (see, e.g., [Z2]). The wonderful new book
by Manuel Kauers and Peter Paule[KP] 
also presents a convincing case.

{\bf Rational Generating Functions}

Equivalently, a $C$-finite sequence is a sequence $\{a(n)\}$ whose 
{\bf ordinary generating function}, $\sum_{n=0}^{\infty} a(n)z^n$,
is rational and where the degree of the denominator is  more
than the degree of the numerator.
These come up a lot in combinatorics and elsewhere
(e.g. formal languages). See the {\it old testament}[St], chapter 4,  and the {\it new testament}[KP], chapter 4.

\vfill\eject

{\bf Why is the Number of tilings of the Knuth Square-Ring C-finite?}

Each of the four ``corners'' ($2 \times 2$ squares) is tiled in  a certain way, where either the tiles
covering it {\it only} cover it, or some tiles also cover neighboring cells {\it not} in the  corner-square.
There are finitely many such scenarios for each corner-square, hence for the
Cartesian product of these scenarios. Having decided how to cover these four corner squares, one must decide
how to tile the remaining four sides, each of which is either a 2 by n rectangle, or one with a few
bites taken from one or both ends. By the {\it transfer matrix method} ([St], 4.7), it follows
{\it a priori} that each of these enumerating sequences is $C$-finite, hence their product,
and hence adding up these finitely many sequences.
It is also easy to see that a rigorous upper bound
for the order is $8$, so the above ``empirical'' approach is justified.

But we still need to be able to compute the first $16$ terms. If we had a very large and very fast
computer, we can actually construct all the tilings, and then count them, but
$104073336269956$ is a pretty big number, so we need a {\it more efficient way}.

{\bf A More Efficient Way}

Suppose that you are given a set of unit-squares (let's call them cells) and you want 
your computer to find the number of ways of tiling it with dominoes. You pick any
cell (for the sake of convenience the left-most bottom-most cell), and look at all the ways
in which to cover it with a domino piece. For each of these ways, removing that tile
leaves a smaller region, thereby getting an obvious  {\it dynamical programming} recurrence.

{\eightrm Procedure {\eighttt NT(R)} in the Maple package {\eighttt RITSUF} computes the number of domino tilings
of any set of cells}.

Using this method, the first-named author solved Knuth's problem in 15 seconds, by typing 
(in a Maple session, in a directory where {\tt RITSUF} has been downloaded to)

{\tt read RITSUF: SeqFrameCsqDirect(2,2,2,2,20,t); } .

We will soon try to solve analogous problems for fatter frames, 
but then we would need to be {\it even} more efficient, and using this more efficient method,
the same calculation would take less than half a second, typing:

{\tt read RITSUF: SeqFrameCsq(2,2,2,2,20,t); } .

{\bf An even More Efficient Way}

Our general problem is to find an efficient automatic way to compute the $C$-finite description, and/or
its generating function (that {\it must} always be a rational function) for the number of
domino tiling of the region, that we denote, in {\tt RITSUF}, by

{\tt Frame(a1,a2,b1,b2,n,n)},

that consists of an $(a1+n+a2) \times (b1+n+b2)$ rectangle  with 
the middle $n \times n$ square removed.

Before describing the algorithm, let us mention that
this is accomplished by procedure

{\tt  SeqFrameC(a1,a2,b1,b2,N,t)} \quad .

For example, {\tt SeqFrameC(1,3,3,1,30,t);} yields:
$$
-4\,{\frac {-9+8\,t+29\,{t}^{2}-10\,{t}^{3}-7\,{t}^{4}+2\,{t}^{5}}{ \left( {t}^{2}-4\,t+1
 \right)  \left( {t}^{4}-4\,{t}^{2}+1 \right) }} \quad .
$$
The fifth argument, $N$ is a parameter for ``guessing'' the $C$-finite description, indicating how many data points to gather
before one tries to guess the $C$-finite description. It is easy to find a priori upper bounds, but it is more fun
to let the user take a guess, and increasing it, if necessary.

{\bf Mihai Ciucu's Amazing Theorem}

The sequence  of positive integers,
demanded by Knuth, enumerating the domino tilings of $Frame(2,2,2,2,n,n)$, 
turned out to be all {\it perfect squares}.
This is not a coincidence! A beautiful theorem  of Mihai Ciucu[C] tells us that whenever
there is a reflective symmetry, the sequence enumerating the dimer counting is either a perfect
square or twice a perfect square. Since we know that (and even if we didn't, we could have
discovered it empirically for each case that we are computing), we have to gather far less
data points. For example, according to the first-named author of this article,
the number of domino tilings of $Frame(3,3,3,3,n,n)$ is $2B(n)^2$, where
$$
\sum_{n=0}^{\infty} B(n)t^n \, = \,
-2\,{\frac {-29+19\,t+102\,{t}^{2}-32\,{t}^{3}-25\,{t}^{4}+7\,{t}^{5}}{ \left( {t}^{2}-4\,t+1 \right)  \left( {t}^{4}-4\,{t}^{2}+1 \right) }} \quad,
$$
while the number of domino tilings of $Frame(4,4,4,4,n,n)$ is $C(n)^2$ where
$$
\sum_{n=0}^{\infty} C(n)t^n= \frac{P(t)}{Q(t)} \quad, 
$$
where
$$
P(t)=
-4 (
901+2517\,t-17574\,{t}^{2}-46322\,{t}^{3}+112903\,{t}^{4}+291045\,{t}^{5}-269376\,{t}^{6}-741508\,{t}^{7}+215233\,{t}^{8}+786069\,{t}^{9}-
$$
$$
21836\,{t}^{10}-352896\,{t}^{11}-24137\,{t}^{12}+67487\,{t}^{13}+5874\,{t}^{14}-5056\,{t}^{15}-359\,{t}^{16}+97\,{t}^{17}) \quad , \quad and \quad
$$ 
$$
Q(t)=
 \left( t-1 \right)  \left( t+
1 \right)  \left( {t}^{4}+{t}^{3}-5\,{t}^{2}+t+1 \right)  \left( {t}^{4}-11\,{t}^{3}+25\,{t}^{2}-11\,t+1 \right)   \cdot
$$
$$
\left( {t}^{4}+7\,{t}^{3}+13\,{t}^{2}+7
\,t+1 \right)  \left( {t}^{4}-{t}^{3}-5\,{t}^{2}-t+1 \right)  \quad .
$$

If you want to see the analogous expressions for $Frame(5,5,5,5,n,n)$ and $Frame(6,6,6,6,n,n)$, then
you are welcome to look at the output file

{\tt http://www.math.rutgers.edu/\~{}zeilberg/tokhniot/oRITSUF2},

that you can generate yourself, {\it ab initio} by running (once you uploaded {\tt RITSUF} onto a Maple session)

{\tt http://www.math.rutgers.edu/\~{}zeilberg/tokhniot/inRITSUF2} \quad .

Our method is {\it pure guessing}, but in order to guess, we need to {\it efficiently} generate sufficiently
many terms of the counting sequence. We must start with {\it rectangles} of fixed width.

{\bf The Number of Domino Tilings of a Rectangle of a Fixed Width}

Let $m$ be a fixed positive integer. We are interested in a
$C$-finite description, as a function of $n$, of the sequence
$A_m(n)$, the number of domino tilings of an $m \times n$
rectangle. In fact, for this {\it specific} problem there
is an ``explicit'' solution, famously found by 
Kasteleyn[Ka] and Fisher \& Temperly [FT], 
but their solution only applies to domino tiling, and we want
to illustrate the {\it general} method.

Also the general approach, using the transfer-matrix method, is
not {\it new} as such (see, e.g. [St]), but since we need it
for counting more elaborate things, let's review it.

Consider the task of tiling the $n$ columns of an $m \times n$ rectangle.
Let's label the cells of a given column from bottom to top
by $\{1, \dots , m \}$.
When we start, all the $m$ cells of the leftmost column are
available, so we start with the {\it state} $\{1, \dots , m\}$.
As we keep on going, not all cells of the current
column are available, since some of them have already been tiled
by the previous column. The other extreme is the empty set,
where nothing is available, and one must go immediately to the
next column, where now everything is available, i.e. the
only follower of $\emptyset$ is the universal set $\{1, \dots , m\}$.

Let's take $m=4$ and see what states may follow the state $\{1,2,3,4\}$.
We may have two vertical tiles  $\{1,2\}$ and $\{3,4\}$, leaving
all the cells of the next column available, yielding the state
$\{1,2,3,4\}$. We may decide instead to {\it only} use horizontal tiles, leaving
nothing available for the next column, resulting in the state $\emptyset$.
If we decide to have one vertical tile in the current column, then 

If it is $\{1,2\}$ then both $3$ and $4$ are parts of horizontal tiles that go
to the next column, leaving the set of cells $\{1,2\}$ available.

If it is $\{2,3\}$ then both $1$ and $4$ are parts of horizontal tiles that go
to the next column, leaving the set of cells $\{2,3\}$ available.

If it is $\{3,4\}$ then both $1$ and $2$ are parts of horizontal tiles that go
to the next column, leaving the set of cells $\{3,4\}$ available.

It follows that the ``followers'' of the state $\{1,2,3,4\}$ are the five states
$$
\emptyset, \{1,2,3,4\}, \{1,2\}, \{2,3\},\{3,4\} \quad .
$$

Who can follow the state $\{1,4\}$? Since only cell $1$ and cell $4$ are available
there can't be any vertical tiles, and both must be parts of horizontal tiles,
occupying $\{1,4\}$ of the next column, and leaving $\{2,3\}$ available,
so the state $\{1,4\}$ has only one follower, the state $\{2,3\}$.

{\eightrm Check out procedure {\eighttt Followers(S,m)} in {\eighttt Ritsfuf}} \quad .

This way we can view any tiling of an $m \times n$ rectangle
as a walk of length $n$ in a directed graph whose vertices are labeled by
subsets of $\{1, \dots , m \}$. The transition matrix of this graph is
the $2^m \times 2^m$ matrix whose rows and columns  naturally correspond to subsets
(in {\tt RITSUF} we made the natural convention that the indices of the rows and columns
correspond to the binary representations of the subsets plus 1).
Let's call $SN(i)$ the set of positive integers corresponding to the positive integer $i$.
For example, $SN(1)$ is the empty set, $SN(10)$  is $\{1,4\}$ etc.

{\eightrm Check out procedure {\eighttt TM(m)} in {\eighttt Ritsfuf} for the transition matrix}.

Calling this transfer matrix $A_m$, the number of domino tilings of an $m \times n$ rectangle is
the $(2^m,2^m)$ entry of the matrix $A_m^n$, since we have to completely tile it,
the starting state must be $\{1, \dots , m\}$, of course, but so is the ending state,
since every thing must be covered, leaving the next column completely available.

{\eightrm Check out procedure {\eighttt SeqRect(m,N)} in {\eighttt Ritsfuf} for the 
counting sequence for the number of domino tilings of an $m$ by $n$ rectangle for $n=0,1, \dots,N$  }.

But not just the $(2^m,2^m)$ entry of $A_m^n$  is informative. Each and every one
of the $(2^m)^2$ entries contains information!
A typical $(i,j)$ entry of $A_m^n$
tells you the number of ways of tiling an $m \times n$ rectangle where the leftmost
column only has the cells in $SN(i)$ available for use, while the rightmost column has
some tiles that stick out, leaving available for the $(n+1)$-th column the
cells of $SN(j)$.

{\bf Counting the Number of Domino Tilings of a Holey Rectangle}

We want to figure out how to use matrix-multiplication to determine the number of tilings of
the region 
$$
Frame(a_1,a_2,b_1,b_2,m,n)
$$
that consists of an $(a_1+m+a_2) \times (b_1+n+b_2)$ rectangle  with 
the middle $m \times n$ rectangle removed.

There are four corner rectangles in our frame:

$\bullet$ the left-bottom (SW) corner consisting of an $a_1 \times b_1$ rectangle \quad ,

$\bullet$ the right-bottom (SE) corner consisting of an $a_1 \times b_2$ rectangle \quad ,

$\bullet$ the left-top (NW) corner consisting of an $a_2 \times b_1$ rectangle \quad ,

$\bullet$ the right-top (NE) corner consisting of an $a_2 \times b_2$ rectangle \quad .

If we look at a typical tiling of the region $Frame(a_1,a_2,b_1,b_2,m,n)$, and focus on
the induced tilings of the four corner-rectangles, we get a tiling with (usually) some
of the tiles intersecting the adjacent non-corner rectangles.

Indeed, for the left-bottom $a_1 \times b_1$ corner-rectangle
(usually) some of the tiles covering its very top row
intersect the very bottom row of the left (West) $m \times b_1$ rectangle,
and  (usually) some of the tiles covering its very right column
intersect the leftmost column of the bottom (South) $a_1 \times n$ rectangle,
and similarly for the other three corner rectangles. So suppose
that the East $m \times b_1$ rectangle has already been tiled,
with  some of its bottom tiles overlapping with our  above-mentioned
left-bottom $a_1 \times b_1$ corner-rectangle, leaving only some
of the cells in the top row available, and
after we complete the tiling of that left-bottom $a_1 \times b_1$ corner-rectangle,
we may use-up some of the  cells of the left column of the bottom (South)
$a_1 \times n$ rectangle, and the complement of that occupied set is only available
for tiling.

This leads naturally to a transfer matrix between the ``states'' of one side of a corner-rectangle to the states
of another side of that corner-rectangle. So let's define $RTM(a,b,S_1,S_2)$, for positive integers $a$ and $b$,
and distinct $S_1,S_2 \in \{1,2,3,4\}$ where we make the convention

1=Up Side \quad , \quad  2=Left Side \quad , \quad  3=Down Side\quad , \quad   4=Right Side \quad ,

that tells you the number of ways of tiling the $a \times b$ rectangle where the
two sides that are not in $\{ S_1, S_2  \}$ are ``smooth'' (i.e. nothing sticks out)
and the two sides $S_1,S_2$ may (and usually do) have some of their tiles ``sticking out''.

Going counterclockwise starting at the left-bottom (SW) corner, we need to find

$\bullet$ For the left-bottom (SW) corner consisting of an $a_1 \times b_1$ rectangle
we need $RTM(a_1,b_1,1,4)$ \quad ,

$\bullet$ For the right-bottom (SE) corner consisting of an $a_1 \times b_2$ rectangle
we need $RTM(a_1,b_2,2,1)$ \quad ,

$\bullet$ For the right-top (NE) corner consisting of an $a_2 \times b_2$ rectangle
we need $RTM(a_2,b_2,3,2)$ \quad ,

$\bullet$ For the left-top (NW) corner consisting of an $a_2 \times b_1$ rectangle
we need $RTM(a_2,b_1,4,3)$ \quad .

Like the matrices $TM(m)$ that for each needed $m$, we only compute {\bf once} 
and then record it, (using {\tt option remember}),
we also only compute $RTM(a,b,S_1,S_2)$ once for each needed $a,b,S_1,S_2$ and remember it.

But how to compute $RTM(a,b,S_1,S_2)$? We first construct, {\it literally}, all the 
domino tilings that completely cover the cells of the $a \times b$ rectangle where
nothing sticks out of the sides that are not labelled $S_1$ or $S_2$, but that
may (and usually do) stick out from the sides labelled $S_1$ and $S_2$. Then we look
at all the  pairs of states, and form a matrix whose $(i,j)$ entry is the
number of stuck-out tilings of the $a \times b$ rectangle where the 
``stick-out'' state of side $S_1$ corresponds to the set labelled $i$, 
and the ``stick-out'' state of side $S_2$ corresponds to the set labelled $j$.

{\eightrm See procedure {\eighttt RTM(a,b,S1,S2)} for its implementation in {\eighttt RITSUF}}.

It follows that, in terms of the matrices $TM(m)$ and $RTM(a,b,S_1,S_2)$, the quantity of
interest, the number of tilings of the  rectangular picture-frame
$Frame(a_1,a_2,b_1,b_2,m,n)$, is
$$
Trace (RTM(a1_1,b_1,1,4)TM(a1)^nRTM(a_1,b_2,2,1)TM(b_2)^mRTm(a_2,b_2,3,2)TM(a_2)^nRTM(a_2,b_1,4,3)TM(b_1)^m ) \quad .
$$

Since matrix power-raising is very fast, and so is matrix-multiplication, we can quickly crank-out sufficiently
many terms in the enumerating sequence, and since we know {\it a priori} that it is $C$-finite, we
can guess a $C$-finite description (and/or a rational generating function), that is proved rigorously
{\it a posteriori} by checking that the order bounds are right.

{\bf The Bivariate Generating Function}

If you are interested in the discrete function of the {\it two} discrete variables $m$ and $n$, 
for the number of domino tilings of  $Frame(a_1,a_2,b_1,b_2,m,n)$, then it is
{\it doubly} $C$-finite, meaning that its bivariate generating function has the form
$P(x,y)/(Q_1(x) Q_2(y))$, for some polynomials $P(x,y),Q_1(x),Q_2(y)$. Using an analogous method for guessing, after cranking-out enough data,
we can get these generating functions easily, using procedure

{\tt GFframeDouble(a1,a2,b1,b2,x,y,N)}

in {\tt RITSUF}. For example, if $D(m,n)$ is the number of domino tilings of
$Frame(2,2,2,2,m,n)$, then
$$
\sum_{m=0}^{\infty}\sum_{n=0}^{\infty} D(m,n)x^my^n \, = \, {\frac{P(x,y)}{Q(x,y)}} \quad ,
$$
where
$$
P(x,y)=
4\,{x}^{3}{y}^{3}-7\,{x}^{3}{y}^{2}-7\,{x}^{2}{y}^{3}-14\,{x}^{3}y+10\,{x}^{2}{y}^{2}-14\,x{y}^{3}+13\,{x}^{3}
$$
$$
+35\,{x}^{2}y+35\,x{y}^{2}+13\,{y}^{3}-30\, {x}^{2}+10\,xy-30\,{y}^{2}-23\,x-23\,y+36 \quad ,
$$
and
$$
Q(x,y)=
\left( x-1 \right)  \left( x+1 \right)  \left( {x}^{2}-3\,x+1 \right)  \left( y-1 \right)  \left( y+1 \right)  \left( {y}^{2}-3\,y+1 \right)  \quad .
$$

{\bf Tiling Crosses}

In how many ways can we tile a cross whose center is a $2 \times 2$ square and each of the four arms have
length $n$? The answer is obtained by typing, in {\tt RITSUF}: {\tt SeqCrossCsq(2,2,20,t);},
getting (in $0.072$ seconds!) that the number is $2B_2(n)^2$, where
$$
\sum_{n=0}^{\infty} B_2(n)t^n \, = \,
{\frac {1}{ \left( t+1 \right)  \left( {t}^{2}-3\,t+1 \right) }} \quad.
$$

Incidentally, $\{B_2(n)\}$  is {\tt http://oeis.org/A001654}, the ``Golden rectangle number''
$F_n F_{n+1}$, so the number of tilings of this cross is $2F_n^2F_{n+1}^2$, and it is
a fairly simple exercise for humans to prove this fact. But we doubt that any human can derive, by hand,
the answer to the analogous question for the cross whose center is a $4 \times 4$ square.
The answer turns out to be $B_4(n)^2$, where
$$
\sum_{n=0}^{\infty} B_4(n)t^n \, = \,
-2\,{\frac {3-t-5\,{t}^{2}+13\,{t}^{3}-11\,{t}^{4}-2\,{t}^{5}+2\,{t}^{6}}{ \left( t-1 \right)  \left( {t}^{4}-11\,{t}^{3}+25\,{t}^{2}-11\,t+1 \right) 
 \left( {t}^{4}+7\,{t}^{3}+13\,{t}^{2}+7\,t+1 \right) }} \quad .
$$
Even SBE needed 4.528 seconds to derive this formula after DZ typed:

{\tt SeqCrossCsq(4,4,30,t);} \quad .

Once again, the approach is {\it purely empirical}. The computer finds all the tilings that completely cover the central square (or rectangle),
possibly (and usually) with some tiles extending beyond. For each such scenario we have the state for each of the
four arms of the cross. Then we use the previously computed (and remembered!) $TM(a)$ matrices, find the corresponding
entry in $TM(a)^n$, and multiply all these four numbers. Finally we (meaning our computers) add up all these
scenarios, getting the desired number. See the source code of {\tt SeqCross(a, b, N)}.

{\bf Monomer-Dimer Tilings}

The beauty of programming is that, once we have finished writing a program, it is easy to 
modify it in order to solve more general, or analogous, problems. By typing {\tt ezraMD();} in {\tt RITSUF}
the readers can find the list of analogous procedures for monomer-dimer tilings, where one can use
either a $1 \times 2$ a $2 \times 1$ or $1 \times 1$ tile, or equivalently, tiling with dominoes
where one is not required to cover all the cells.

For example, if $A(n)$ is the number of ways of tiling the Knuth region 
(obtained by removing the central $n^2$ squares from an $(n+4) \times (n+4)$ array of squares)
by dimers {\it and} monomers, the answer is much messier.
It took SBE $15$ seconds to discover that
$$
\sum_{n=0}^{\infty} A(n)t^n \, = \,
\frac{P(t)}{Q(t)} \quad,
$$
where
$$
P(t)=
-94\,{t}^{30}+1361\,{t}^{29}+43975\,{t}^{28}-494267\,{t}^{27}-5787443\,{t}^{26}+61186056\,{t}^{25}+266911158\,{t}^{24}
$$
$$
-3200500450\,{t}^{23}-3505671568\,
{t}^{22}
+74767156291\,{t}^{21}-29007687275\,{t}^{20}-796609853769\,{t}^{19}+823823428983\,{t}^{18}
$$
$$
+3924729557742\,{t}^{17}-4977782472712\,{t}^{16}-9040256915004\,{t}^{15}+11643454084810\,{t}^{14}+9751493606823\,{t}^{13}
$$
$$
-11693567793807\,{t}^{12}-4837640809485\,{t}^{11}+5123918478955\,{t}^{10}+1059903067708\,{t}^{9}-944330286322\,{t}^{8}-87120095554\,{t}^{7}
$$
$$
+67451928324\,{t}^{6}+657867045\,{t}^{5}-1679236205\,{t}^{4}+73176689\,{t}^{3}+6962033\,
{t}^{2}-226706\,t-10012 \quad ,
$$
and
$$
Q(t)=
 \left( t-1 \right)  \left( {t}^{3}-7\,{t}^{2}+11\,t-1 \right)  \left( {t}^{3}+7\,{t}^{2}-33\,t-1 \right)  \left( {t}^{3}-11\,{t}^{2}+7\,t-1 \right)  \cdot
$$
$$
 \left( {t}^{3}+{t}^{2}-3\,t-1 \right)  \left( {t}^{3}-27\,{t}^{2}+107\,t-1 \right)  \left( {t}^{3}+3\,{t}^{2}-t-1 \right)  \cdot
$$
$$
\left( {t}^{6}+20\,{t}^{5}+55
\,{t}^{4}-304\,{t}^{3}-337\,{t}^{2}+8\,t+1 \right)  \left( {t}^{6}-37\,{t}^{4}-76\,{t}^{3}-37\,{t}^{2}+1 \right) \quad .
$$

\vfill\eject

{\bf Weighted Counting}

What if instead of {\it straight counting} one wants to do {\it weighted counting}?.
The weight of a domino tiling
is defined to be  $h^{\#HorizontalTiles}v^{\#VerticalTiles}$ (in the case of domino tilings) and \hfill\break
$h^{\#HorizontalTiles}v^{\#VerticalTiles}m^{\#MonomerTiles}$ (in the case of monomer-dimer tiling).

For this we have the Maple package {\tt RITSUFwt}, freely downloadable from \hfill\break
{\eighttt http://www.math.rutgers.edu/\~{}zeilberg/tokhniot/RITSUFwt}, 
that also contains all the procedures of {\tt RITSUF}. To get a list of the procedures
for weighted counting, type:

{\tt ezraWt(); }

For example, to get the weighted analog of the Knuth problem, type:

{\tt SeqFrameCwt(2,2,2,2,30,t,h,v);}

To see the output, go to:

{\tt http://www.math.rutgers.edu/\~{}zeilberg/tokhniot/oRITSUFwt1} \quad ,

where one can also see statistical analysis.

{\bf The Maple package} {\tt ARGF}

The Maple package {\tt ARGF} (short for ``Analysis of Rational Generating Functions'') downloadable from:

{\eighttt http://www.math.rutgers.edu/\~{}zeilberg/tokhniot/ARGF} \quad ,

does automatic statistical analysis of random variables whose weight-enumerators
are given by rational  functions, like the one outputted by {\tt RITSUFwt}, and
whose procedures are also included in the latter. We urge the readers to look at

{\tt http://www.math.rutgers.edu/\~{}zeilberg/tokhniot/oRITSUFwt4} \quad ,

and the other output files in the webpage of this article

{\tt http://www.math.rutgers.edu/\~{}zeilberg/mamarim/mamarimhtml/ritsuf.html} \quad ,

for examples. Here we use the methodology of [Z3][Z4].

{\bf Conclusion}

The {\it deductive method} in mathematics reigned for the last 2500 years. It is time to
replace it by the {\it inductive method}. Often the inductive method is fully rigorous
(like here), but other times, it may not be, but {\it who cares}?. Only trivial results
can be proved fully rigorously, so let's not tie ourselves with these antiquated shackles,
and explore mathematics experimentally!

{\bf References}

[C] Mihai Ciucu, {\it Enumeration of perfect matchings in graphs with reflective symmetry},
 J. Combin. Theory Ser. A {\bf 77}(1997), 67-97.

[FT] M. Fisher and H. Temperley, {\it Dimer Problems in Statistical
Mechanics-an exact result}, Philos. Mag. {\bf 6} (1961), 1061-1063.

[Ka] P. W. Kasteleyn, {\it The statistics of dimers on a lattice: I. The number of
dimer arrangements in a quadratic lattice}, Physica {\bf 27} (1961), 1209-1225.

[KP] Manuel Kauers  and Peter Paule, {\it ``The Concrete Tetrahedron''}, Springer, 2011.

[Kn] Donald E. Knuth, {\it Tiling a square ring with dominoes}, proposed problem
\#1868, Math. Magazine {\bf 84 No. 2} (April 2011) Solution by
John Bonomo and David Offner,  Math. Magazine {\bf 85 No. 2} (April 2012), 154-155.

[St] Richard Stanley, { \it ``Enumerative combinatorics''}, volume {\bf 1}, Wadsworth and Brooks/Cole, Pacific Grove, CA, 1986,
second printing, Cambridge University Press, Cambridge, 1996.

[T] Roberto Tauraso, {\it A New Domino Tiling Sequence}, Journal of Integer Sequences, {\bf 7} (2004), Article 04.2.3 .Available on-line: \hfill\break
{\tt http://www.cs.uwaterloo.ca/journals/JIS/VOL7/Tauraso/tauraso3.pdf} .

[Z1] Doron Zeilberger, {\it  The C-finite Ansatz}, a ``sanitized'' version is to appear in the Ramanujan Journal, the
original, uncensored version is available on-line: \hfill\break
{\tt http://www.math.rutgers.edu/\~{}zeilberg/mamarim/mamarimhtml/cfinite.html}

[Z2] Doron Zeilberger, 
{\it  An Enquiry Concerning Human (and Computer!) [Mathematical] Understanding},
Appeared in: C.S. Calude, ed., ``Randomness \& Complexity, from Leibniz to Chaitin'' World Scientific, Singapore, 2007.
Available on-line: \hfill\break
{\tt http://www.math.rutgers.edu/\~{}zeilberg/mamarim/mamarimhtml/enquiry.html}

[Z3] Doron Zeilberger, 
{\it The Automatic Central Limit Theorems Generator (and Much More!)},
in: 
``{\it Advances in Combinatorial Mathematics}'',
Proceedings of the Waterloo Workshop in Computer Algebra 2008 in honor of Georgy P. Egorychev, 
chapter 8, pp. 165-174, (I.Kotsireas, E.Zima, eds.), Springer Verlag, 2009.
Available on-line: \hfill\break
{\tt http://www.math.rutgers.edu/\~{}zeilberg/mamarim/mamarimhtml/georgy.html}

[Z4] Doron Zeilberger, 
{\it HISTABRUT: A Maple Package for Symbol-Crunching in Probability theory}, Personal Journal of Shalosh B. Ekhad
and Doron Zeilberger, Aug. 25, 2010, \hfill\break
{\tt http://www.math.rutgers.edu/\~{}zeilberg/mamarim/mamarimhtml/histabrut.html}
\end